\documentclass[12pt,a4paper]{article}
\usepackage{amsmath,amsfonts,amssymb,amsthm}
\usepackage{mathrsfs,dsfont}

\newtheorem{defi}{Definition}[section]

\newtheorem{prop}[defi]{Proposition}
\newtheorem{thm}[defi]{Theorem}
\newtheorem{cor}[defi]{Corollary}
\theoremstyle{definition}
\newtheorem{rem}[defi]{Remark}

\def\C{\mathds C}
\def\phi{\varphi}
\def\rho{\varrho}
\def\id{{\it id}}
\def\ul{\underline}
\def\A{\mathds A}
\def\O{\mathcal O}

\DeclareMathOperator{\KG}{\textit{KG}}
\DeclareMathOperator{\im}{im}
\DeclareMathOperator{\Spec}{Spec}
\DeclareMathOperator{\GL}{GL}
\DeclareMathOperator{\SL}{SL}
\DeclareMathOperator{\Mor}{Mor}

\newcommand{\catXschemes}[1]{(\textit{$#1$-schemes})}
\def\catsets{(\textit{sets})}

\DeclareMathOperator{\FGrass}{\underline{Grass}}
\DeclareMathOperator{\Grass}{Grass}
\DeclareMathOperator{\FGGrass}{\underline{G-Grass}}
\DeclareMathOperator{\GGrass}{G-Grass}
\DeclareMathOperator{\Hilb}{Hilb}
\DeclareMathOperator{\FGHilb}{\underline{G-Hilb}}
\DeclareMathOperator{\GHilb}{G-Hilb}

\begin{document}

\title{Construction of $G$-Hilbert schemes}
\author{Mark Blume}
\date{}

\maketitle

\begin{abstract}
In this paper we construct $G$-Hilbert schemes for finite
group schemes $G$.
We find a construction of $G$-Hilbert schemes as relative $G$-Hilbert 
schemes over the quotient that does not need the Hilbert scheme of $n$ 
points, works under more natural assumptions and gives additional 
information about the morphism from the $G$-Hilbert scheme to the quotient.
\end{abstract}
                
\bigskip

\section*{Introduction}

For a quasi-projective variety $X$ over $\C$ and $G$ a finite group
of automorphisms, the $G$-Hilbert scheme $\GHilb_\C X$, a variation 
of the Quot and Hilbert scheme construction \cite{Gr61}, has been 
defined as the scheme that parametrises $G$-stable finite closed 
subschemes $Z\subseteq X$ with coordinate ring isomorphic to the 
regular representation.
These subschemes $Z$ are sometimes called $G$-clusters. A free 
$G$-orbit is a $G$-cluster, any $G$-cluster is supported by a 
$G$-orbit, but in the case of nontrivial stabiliser there might 
be many nonreduced $G$-clusters supported by the same orbit.
The $G$-Hilbert scheme is related to the quotient by a morphism 
$\tau\colon\GHilb_\C X\to X/G$ taking $G$-clusters to the 
corresponding orbits. This morphism is projective, and over the 
set of free $G$-orbits it is an isomorphism. 

The $G$-Hilbert scheme has been introduced to resolve quotient 
singularities and to describe the properties of such resolutions:
for finite subgroups $G\subset\GL(2,\C)$ the morphism 
$\tau\colon\GHilb_\C X\to X/G$ is the minimal resolution, in 
some higher dimensional cases including finite subgroups 
$G\subset\SL(3,\C)$ a crepant resolution \cite{BKR01}.
This has been used extensively in works concerning the
McKay correspondence (see \cite{Re97}, \cite{Re99} for surveys), 
in particular it is essential in the proof of the derived 
McKay correspondence in \cite{BKR01}. 

In the original definition \cite{ItNm96}, \cite{ItNm99} the 
Hilbert scheme of $G$-orbits is defined as closure of the 
set of free $G$-orbits in the fixed point subscheme of the 
Hilbert scheme of $|G|$ points.
In this paper we use the slightly different definition as a moduli 
space of $G$-clusters, i.e.\ the $G$-Hilbert scheme is the scheme 
that represents the corresponding moduli functor. This first appears 
in \cite{Re97} and is adopted in many of the following works.
To compare the two definitions, the Hilbert scheme of $G$-orbits as 
in the original definition is the component of the moduli space of 
$G$-clusters birational to the quotient, see also 
\cite[section 4.1]{CrRe02}.
$G$-Hilbert schemes for abelian groups $G$ are further studied 
in \cite{CMT07a}, \cite{CMT07b}, their results include an example 
of a reducible $G$-Hilbert scheme.

Usually the $G$-Hilbert scheme and the morphism to the quotient 
$\GHilb X\to X/G$ are constructed using the Hilbert scheme
of $n$ points and the morphism\linebreak $\Hilb^n X\to S^nX$
to the symmetric product, which is part of the general framework 
of Hilbert-to-Chow morphisms.
But it is of interest to give a construction of the $G$-Hilbert 
scheme independent of the Hilbert scheme of $n$ points.
$G$-Hilbert schemes also appear as moduli spaces of representations
of quivers and as such can be constructed as GIT quotients, see 
e.g.\ \cite{CMT07a} for the case of abelian groups $G$ acting 
linearly on affine spaces.

The new idea carried out in the present paper is to construct 
the $G$-Hilbert scheme as a scheme over the quotient. This is
done by embedding it into some $G$-equivariant Grassmannian.
A construction close to the construction in the present paper 
but under more special assumptions appeared later in \cite{Ga07}. 

As one result of this paper we give a new construction of the
$G$-Hilbert scheme and the morphism to the quotient that does
not need the Hilbert scheme of $n$ points (theorem \ref{thm:reprGHilb}). 
We introduce relative $G$-Hilbert schemes and we construct the 
$G$-Hilbert scheme $\GHilb_KX$ as a relative $G$-Hilbert scheme 
$\GHilb_{X/G}X$ over the quotient, this way we avoid the embedding 
into the Hilbert scheme of $n$ points. 
It is done under conditions more general and more natural as before, 
in particular, no quasi-projectivity assumption will be needed during 
the construction. Further, we obtain additional information about the 
natural morphism $\tau\colon\GHilb_KX\to X/G$, which is identified
with the structure morphism of the relative $G$-Hilbert scheme 
$\GHilb_{X/G}X$: we have the commutative diagram 
(remark \ref{rem:GHilboverX/G})

\medskip

\noindent
\begin{picture}(150,26)
\put(60,20){\makebox(0,0)[r]{$\GHilb_{X/G}X$}}
\put(75,21){\makebox(0,0)[b]{$\sim$}}
\put(90,20){\makebox(0,0)[l]{$\GHilb_KX$}}
\put(75,20){\vector(1,0){10}}\put(75,20){\vector(-1,0){10}}
\put(90,16){\vector(-1,-1){8}}\put(60,16){\vector(1,-1){8}}
\put(87,10){\makebox(0,0)[l]{\small$\tau$}}
\put(75,2){\makebox(0,0)[b]{$X/G$}}
\end{picture}

\medskip

Using the construction of the $G$-Hilbert scheme as a relative 
$G$-Hilbert scheme it is easy to see that the fibres of $\tau$
are again $G$-Hilbert schemes (remark \ref{rem:fibresGHilb->quot}),
and it is possible to apply the differential study of Quot 
schemes \cite{Gr61} to calculate relative tangent spaces
of $\GHilb_K\A^n_K$ over $\A^n_K/G$. As argued in \cite{Bl07},
there are relations to the stratification of the $G$-Hilbert 
scheme considered in \cite{ItNm99}. 

\pagebreak

In another direction we generalise the $G$-Hilbert scheme 
construction to finite group schemes over fields not necessarily 
algebraically closed or of positive characteristic
(result again formulated in theorem \ref{thm:reprGHilb}). 
This generalisation is motivated by generalisations of the McKay 
correspondence \cite{McK80} to base schemes other than $\Spec\C$;
the case of non algebraically closed fields $K$ of characteristic 
$0$ and finite subgroup schemes $G\subset\SL(2,K)$ has been 
investigated in \cite{Bl06}. 

\bigskip

\noindent This paper is organised as follows.

In section 1 we define the $G$-Hilbert scheme $\GHilb_SX$ in terms 
of its functor in the generality used in this paper.
We construct the morphism $\tau\colon\GHilb_S X\to X/G$ as a morphism 
of functors (proposition \ref{prop:GHilbX->Y}).

In section 2 we study ways to change the base scheme of 
$G$-Hilbert schemes. The results give some freedom in choosing 
the base scheme. The maximal possible base scheme, i.e.\ with 
minimal fibres, for the $G$-Hilbert functor of a given 
$G$-scheme $X$ is the quotient $X/G$. 
The main result of this section (theorem \ref{thm:GHilb_S-GHilb_Y}) 
in particular implies that the natural morphism 
$\tau\colon\GHilb_S X\to X/G$ can be identified with the 
structure morphism of a relative $G$-Hilbert scheme with the 
quotient as base scheme (remark \ref{rem:GHilboverX/G}).

In section 3 we introduce equivariant and $G$-Grassmannian
functors and show representability of these functors 
(theorem \ref{thm:reprGGrass}).
The other main result (theorem \ref{thm:GHilb->GGrass}) shows, 
that in certain cases there is a natural closed embedding 
of $G$-Hilbert functors into $G$-Grassmannian functors.
Together, these results imply that $G$-Hilbert functors
$\FGHilb_SX$ for finite $X\to S$ are represented by projective 
$S$-schemes (corollary \ref{cor:reprGHilb,X->Sfinite}).

Section 4 contains the theorem about representability of the 
$G$-Hilbert functor (theorem \ref{thm:reprGHilb}). 
Combining the results of sections 2 and 3, we give under natural 
assumptions a construction of the $G$-Hilbert scheme $\GHilb_KX$ 
for algebraic $K$-schemes $X$ as relative $G$-Hilbert scheme 
over the quotient.

The appendix develops some facts about $G$-equivariant sheaves and 
decomposition into isotypic components, which are essential in order 
to extend the $G$-Hilbert scheme construction to group schemes.

\bigskip

{\it Notations, conventions and general definitions.}
In the following let $G=\Spec A$ be a finite group scheme over 
a field $K$ with $p\colon G\to\Spec K$ the structure morphism and
$e\colon\Spec K\to G$ the unit, set $|G|=\dim_KA$. 
A $G$-scheme over a $K$-scheme $S$ is an $S$-scheme with an 
operation of $G$ over $S$.  
In general we write a lower index for base extensions, e.g.\ 
if $X,T$ are $S$-schemes then $X_T$ denotes the $T$-scheme 
$X\times_ST$. For an $\O_S$-module $\mathscr F$ we write 
$\mathscr F_T$ for its pull-back with respect to $T\to S$,
more generally for an $\O_X$-module $\mathscr F$ write
$\mathscr F_T$ for its pull-back with respect to $X_T\to X$.
Functors $\catXschemes{S}^\circ\to\catsets$ we denote as 
$S$-functors.
Sometimes we will underline a functor to distinguish it from 
the corresponding scheme.

\pagebreak
\section{{\boldmath$G$}-Hilbert functors and the morphism to 
the quotient}\label{sec:GHilb}

For a $\,G$-scheme $X$ over $S$ we define the $G$-Hilbert 
functor $\FGHilb_S X\colon\catXschemes{S}^\circ$ $\to\catsets$ by
\[
\FGHilb_SX(T):=
\left\{
\begin{array}{l}
\textit{Quotient $G$-sheaves}\;\;[0\to\mathscr I\to\O_{X_T}\to
\O_Z\to 0]\;\;\textit{on $X_T$},\\
Z\;\textit{finite flat over $T$},\;\;
\textit{for $t\in T\colon$ $H^0(Z_t,\O_{Z_t})$ isomorphic}\\ 
\textit{to the regular representation}\\
\end{array}
\right\}
\]
Here a quotient $G$-sheaf $[0\to\mathscr I\to\O_{X_T}\to\O_Z\to 0]$ 
is an exact sequence $0\to\mathscr I\to\O_{X_T}\to\O_Z\to 0$ 
of quasi-coherent $G$-sheaves on $X_T$ with $\mathscr I,\O_Z$ 
specified up to isomorphism, that is either a quasi-coherent 
$G$-subsheaf $\mathscr I\subseteq\O_{X_T}$ or an equivalence 
class $[\O_{X_T}\to\O_Z]$ of surjective equivariant homomorphisms 
of quasi-coherent $G$-sheaves with two of them equivalent if their 
kernels coincide.
Sometimes we shortly write $Z\in\FGHilb_SX(T)$. A $G$-cluster is 
an element $Z\in\FGHilb_SX(\Spec L)$ for a morphism of $S$-schemes 
$\Spec L\to X$, where $L$ is an extension field of $K$.

For $X$ quasi-projective over $S=\Spec\C$ and $G$ a finite group this 
coincides with the usual definition of the $G$-Hilbert scheme as  
moduli space of $G$-clusters used in the literature, e.g. \cite{Re97} 
(making for $Z\to T$ use of the fact that proper and quasi-finite 
implies finite, \cite[III (1), (4.4.2)]{EGA} or 
\cite[IV (3), (8.11.1)]{EGA}).

There is the natural morphism $\GHilb_SX\to X/G$ taking $G$-clusters 
to the corresponding orbits. It is implied by the following result,
here we work with the functors.

\begin{prop}\label{prop:GHilbX->Y}
Let $Y$ be a $G$-scheme over $S$ with trivial $G$-operation
and\linebreak $\phi\colon X\to Y$ an equivariant morphism of 
$G$-schemes over $S$. Then there is a unique morphism of $S$-functors
\[\tau\colon\FGHilb_S X\to\ul{Y}\]
such that for $S$-schemes $T$ and $Z\in\FGHilb_SX(T)$
with image $\tau(Z)\in Y(T)$ the following diagram
commutes:

\noindent
\begin{picture}(150,24)
\put(65,20){\makebox(0,0)[c]{$Z$}}
\put(85,20){\makebox(0,0)[c]{$X$}}
\put(65,4){\makebox(0,0)[c]{$T$}}
\put(85,4){\makebox(0,0)[c]{$Y$}}
\put(69,20){\vector(1,0){12}}\put(69,4){\vector(1,0){12}}
\put(65,16){\vector(0,-1){8}}\put(85,16){\vector(0,-1){8}}
\put(75,5){\makebox(0,0)[b]{\footnotesize$\tau(Z)$}}
\end{picture}
\vspace{-18mm}
\begin{equation}\label{eq:GHilbX->Y}\,\end{equation}
\vspace{4mm}
\end{prop}
\begin{proof}
The existence and uniqueness of a morphism $\tau(Z)$ such that 
diagram (\ref{eq:GHilbX->Y}) commutes derives from the fact, 
that any $Z\to T$ for $Z\in\FGHilb_SX(T)$ is a geometric and 
thus as well categorical quotient of $Z$ by $G$.
The maps $\FGHilb_SX(T)\to Y(T)$, $Z\mapsto\tau(Z)$ are functorial 
in $T$, that is $\tau(Z_{T'})=\tau(Z)\circ\alpha$ for 
morphisms $\alpha\colon T'\to T$ of $S$-schemes (by uniqueness of 
$\tau(Z_{T'})$).
\end{proof}

\section{Varying the base scheme of {\boldmath$G$}-Hilbert functors}
\label{sec:varybasescheme}

We fix notations and state some generalities on changing the base 
scheme of functors with respect to a morphism $\phi\colon S'\to S$.

Base restriction: for an $S'$-functor $F'$ we define the $S$-functor 
$_SF'$ by taking disjoint unions
\vspace{-2mm}
\[_SF'(T):=\bigsqcup\big\{F'(\alpha)\:\big|\:\alpha\in\Mor_S(T,S')\big\}\]
For an $S'$-scheme $Y$ the corresponding construction is 
to consider $Y$ as an $S$-scheme by composing the structure 
morphism $Y\to S'$ with $S'\to S$.

Base extension: for an $S$-functor $F$ we define the $S'$-functor 
$F_{S'}$ by\linebreak 
$F_{S'}(\alpha\colon\!T\!\to\!S')\!=\!F(\phi\circ\alpha)$,
which is the restriction of $F$ to the category of $S'$-schemes.
$F_{S'}$ also can be realised as the fibred product
$F\times_{\ul{S}}\ul{S'}$ considered as $S'$-functor. 
For an $S$-scheme $X$ this corresponds to the usual base extension 
$X_{S'}=X\times_SS'$.

If in addition a morphism of $S$-functors $\psi\colon F\to\ul{S'}$ 
is given, then we can consider $F$ as a functor on the category 
of $S'$-schemes via $\psi$: let the 
$S'$-functor $F_{(S',\psi)}$ be given by 
\vspace{-2mm}
\[\;F_{(S',\psi)}(\alpha\colon T\to S')
=\{\beta\in F(\phi\circ\alpha)\:|\:\psi(\beta)=\alpha\}
=\{\beta\in F_{S'}(\alpha)\:|\:\psi(\beta)=\alpha\}\]
For schemes this means to consider an $S$-scheme $X$ 
as an $S'$-scheme via a given $S$-morphism $X\to S'$.
Note that $F_{(S',\psi)}$ can be considered as subfunctor 
of $F_{S'}$ and that $_S(F_{(S',\psi)})\cong F$.

\begin{rem}\label{rem:baseextGHilb} (Base extension).
Let $X$ be a $G$-scheme over $S$ and $S'$ an $S$-scheme.
Then there is the isomorphism of $S'$-functors
$(\FGHilb_SX)_{S'}\cong\FGHilb_{S'}X_{S'}$,
which derives from the natural isomorphisms 
$X\times_S T\cong X_{S'}\times_{S'}T$ for $S'$-schemes $T$.
\end{rem}

\begin{thm}\label{thm:GHilb_S-GHilb_Y}
Let $Y$ be a $G$-scheme over $S$ with trivial $G$-operation and 
$\phi\colon X\to Y$ an equivariant morphism of $G$-schemes over $S$. 
Let $\tau\colon\FGHilb_SX\to\ul{Y}$ be the morphism of proposition
\ref{prop:GHilbX->Y}. 
Then there is an isomorphism of $Y$-functors
\[(\FGHilb_SX)_{(Y,\tau)}\,\cong\;\FGHilb_YX\]
\end{thm}
\begin{proof}
Both functors are subfunctors of $(\FGHilb_SX)_Y\colon
\catXschemes{Y}^\circ\to\catsets$:
\[
\begin{array}{l}
(\FGHilb_SX)_{(Y,\tau)}(\alpha\colon T\to Y)
=\{Z\in(\FGHilb_SX)_{Y}(\alpha\colon T\to Y)\:|\:
\tau(Z)=\alpha\}\\[0.4em]
\FGHilb_{Y}X(T)
=\{Z\in(\FGHilb_SX)_Y(T)\:|\:\textit{$Z\hookrightarrow X\times_S T$
factors through $X\times_YT$}\}\\
\end{array}
\]
We show that they coincide. Let $T$ be a $Y$-scheme and 
$Z\subseteq X\times_S T$ a closed subscheme defining an element 
of $(\FGHilb_SX)_{Y}(T)$. Then there are the equivalences
\[
\begin{array}{rcl}
Z\in\FGHilb_YX(T)&\Longleftrightarrow&
\textit{$Z\hookrightarrow X\times_S T$ factors through $X\times_YT$}\\
&\hspace{-4cm}\Longleftrightarrow&
\hspace{-2cm}\textit{diagram {\rm(\ref{eq:GHilbX->Y})} 
commutes for $Z$ and the given morphism $T\to Y$}\\
&\hspace{-4cm}\Longleftrightarrow&
\hspace{-2cm}Z\in(\FGHilb_SX)_{(Y,\tau)}(T)\\
\end{array}
\]
For the last equivalence one uses uniqueness of a morphism $T\to Y$
making diagram (\ref{eq:GHilbX->Y}) commute.
\end{proof}

\begin{rem}\label{rem:GHilboverX/G}
The quotient morphism $X\to X/G$ induces $\tau\colon\FGHilb_SX
\to\ul{X/G}$ (proposition \ref{prop:GHilbX->Y}) and there is 
the isomorphism $(\FGHilb_SX)_{(X/G,\tau)}\cong\FGHilb_{X/G}X$ 
(theorem \ref{thm:GHilb_S-GHilb_Y}). This implies an isomorphism
of $S$-functors $\FGHilb_SX\cong\linebreak\:_S(\FGHilb_{X/G}X)$.
It follows, that if $\FGHilb_{X/G}X$ is representable, 
then so is\linebreak $\FGHilb_SX$ and in this case the 
corresponding schemes $\GHilb_{X/G}X$ and $\GHilb_SX$ are 
isomorphic as $X/G$-schemes ($\GHilb_SX$ as $X/G$-scheme via $\tau$).
This implies an isomorphism of $S$-schemes
$\GHilb_SX\cong{}_S(\GHilb_{X/G}X)$.
\end{rem}

\begin{rem}\label{rem:fibresGHilb->quot}
For points $y\in X/G$ one can determine the fibres of the
natural morphism $\tau\colon\GHilb_SX\to X/G$ (assuming 
representability).
By remark \ref{rem:GHilboverX/G} the fibres of $\tau$ are the fibres 
$(\GHilb_{X/G}X)_y$ of the $X/G$-scheme $\GHilb_{X/G}X$, 
by remark \ref{rem:baseextGHilb} these are isomorphic to 
$G$-Hilbert schemes $\GHilb_{\kappa(y)}X_y$ over $\kappa(y)$,
where the $G$-schemes $X_y$ are the fibres of $X\to X/G$ over $y$.
\end{rem}

\begin{rem} 
Variation of the base scheme can be useful to derive the 
following result (contained in \cite{Te04}) in a simple way. 
For $X,Y$ $G$-schemes over $S$, $Y$ with trivial $G$-operation,
there is the isomorphism 
\[\FGHilb_SX\times_S Y\cong(\FGHilb_SX)\times_{\ul{S}}\ul{Y}\]
Proof. The projection $X\times_S Y\to Y$ is $G$-equivariant, theorem 
\ref{prop:GHilbX->Y} then constructs a morphism of $S$-functors 
$\tau\colon\FGHilb_SX\times_S Y\to\ul{Y}$, by theorem 
\ref{thm:GHilb_S-GHilb_Y} there is an isomorphism of $Y$-functors 
$(\FGHilb_SX\times_S Y)_{(Y,\tau)}\cong\FGHilb_YX\times_S Y$ 
and $\FGHilb_YX\times_S Y\cong(\FGHilb_SX)\times_{\ul{S}}\ul{Y}$
by the usual base extension. Restricting the base to $S$ yields 
the result.
\end{rem}

\section{{\boldmath$G$}-Hilbert schemes and {\boldmath$G$}-Grassmannians}
\label{sec:GHilbGGrass}

The Grassmannian functors $\FGrass_S^n(\mathscr F)$ for a 
quasi-coherent sheaf $\mathscr F$ on a scheme $S$ as defined in 
\cite[\S 9.7]{EGA1} are the $S$-functors given by 
\[\FGrass_S^n(\mathscr F)(T)=
\left\{
\begin{array}{l}
\textit{Quotient sheaves}\;\;[0\to\mathscr H\to\mathscr F_T\to
\mathscr G\to 0]\;\;\textit{on T},\\
\mathscr G\;\textit{locally free of rank $n$}\\
\end{array}
\right\}\]
Similarly define the equivariant Grassmannian functors 
$\FGrass_S^{G,n}(\mathscr F)$ for a quasi-coherent $G$-sheaf 
$\mathscr F$ on $S$ by
\[
\FGrass_S^{G,n}(\mathscr F)(T)=
\left\{
\begin{array}{l}
\textit{Quotient $G$-sheaves}\;\;[0\to\mathscr H\to\mathscr F_T\to
\mathscr G\to 0]\;\;\textit{on T},\\
\mathscr G\;\textit{locally free of rank $n$}\\
\end{array}
\right\}
\]
and the $G$-Grassmannian functor $\FGGrass_S(\mathscr F)$ with 
condition ``$\,\mathscr G$ {\it locally free of finite rank with 
fibres isomorphic to the regular representation}'' replacing 
``$\,\mathscr G$ {\it locally free of rank $n$}".

\begin{rem}
For a $G$-sheaf $\mathscr F$ there is a natural $G$-operation 
on $\FGrass_S^n(\mathscr F)$ given by 
$G(T)\times_S\FGrass_S^n(\mathscr F)(T)\to\FGrass_S^n(\mathscr F)(T)$,
$(g,[0\to\mathscr H\to\mathscr F_T\to\mathscr G\to 0])\mapsto
[0\to g_*\mathscr H\to\mathscr F_T\to g_*\mathscr G\to 0]$
by applying $g_*$ for the automorphism $g\colon X_T\to X_T$ and 
using the isomorphism $\mathscr F_T\to g_*\mathscr F_T$ coming 
from the $G$-sheaf structure of $\mathscr F$ (see appendix 
\ref{sec:Gsheaves}).
\end{rem}

\begin{thm}\label{thm:reprGGrass}
Assume that the Hopf algebra $A$ of $G$ is cosemisimple. 
Then for a quasi-coherent $G$-sheaf $\mathscr F$ on $S$ 
the $G$-Grassmannian functor $\FGGrass_S(\mathscr F)$ 
is represented by an $S$-scheme $\GGrass_S(\mathscr F)$.
If moreover $\mathscr F$ is of finite type, then 
$\GGrass_S(\mathscr F)$ is projective over $S$.
\end{thm}
\begin{proof}
The Grassmannian functor $\FGrass_S^{|G|}(\mathscr F)$ is 
representable \cite[I, (9.7.4)]{EGA1}.

The equivariant Grassmannian $\FGrass_S^{G,|G|}(\mathscr F)$ is a 
closed subfunctor of\linebreak $\FGrass_S^{|G|}(\mathscr F)$, because it
coincides with the fixed point subfunctor \cite[II, \S1, Def. 3.4]{DG} 
with respect to the natural $G$-operation on 
$\FGrass_S^{|G|}(\mathscr F)$ (a $T$-valued point of the fixed point
subfunctor, i.e. a $G$-stable subsheaf $\mathscr H\subseteq\mathscr F_T$,
is the same as a $T$-valued point of the equivariant Grassmannian
functor, i.e.  an injective $G$-equivariant homomorphism $\mathscr H
\to\mathscr F_T$, by remark \ref{rem:Gsubsheaf}). The fixed point 
subfunctor is closed by \cite[II, \S1, Thm. 3.6]{DG}.
 
The $G$-Grassmannian is a component of the equivariant Grassmannian 
\linebreak
$\FGrass_S^{G,|G|}(\mathscr F)$: specifying in addition the 
isomorphism class of a representation determines an open and closed 
subfunctor, since the quotient sheaf $\mathscr G$ has a decomposition 
into a direct sum of isotypic components and the ranks of these
locally free sheaves determine the isomorphism class of the 
representations on the fibres (remark \ref{rem:isotypdecomp}).

If $\mathscr F$ is quasi-coherent of finite type, then the
Grassmannian $\Grass_S^{|G|}(\mathscr F)$ is projective (in the
sense of \cite[II, (5.5)]{EGA}) via the Pl\"ucker embedding 
\cite[I, (9.8.4)]{EGA1}, and this implies that the $G$-Grassmannian 
$\GGrass_S(\mathscr F)$ is projective over $S$.
\end{proof}

Let $X$ be a $G$-scheme over $S$. Assume that $X$ is affine over 
$S$, then $X=\Spec_S\mathscr B$ for a quasi-coherent $G$-sheaf of 
$\O_S$-algebras.
One may rewrite the $G$-Hilbert functor for $X$ over $S$ as  
\[
\FGHilb_SX(T)=
\left\{
\begin{array}{l}
\textit{Quotient $G$-sheaves}\;\;[0\to\mathscr I\to
\mathscr B_T\to\mathscr C\to 0]\\
\textit{of $\mathscr B_T$-modules on $T$},\;\;
\textit{$\mathscr C$ locally free of finite rank}\\
\textit{with fibres isomorphic to the regular representation}\\
\end{array}
\right\}
\]
There is a morphism of $S$-functors $\FGHilb_SX\to\FGGrass_S(\mathscr B)$ 
consisting of injective maps
\[\FGHilb_SX(T)\to\FGGrass_S(\mathscr B)(T),\quad
[\mathscr B_T\to\mathscr C]\mapsto[\mathscr B_T\to\mathscr C]\]  
defined by forgetting the algebra structure of $\mathscr B$,
this way $\FGHilb_SX$ becomes a subfunctor of $\FGGrass_S(\mathscr B)$.
The essential point we will show is the closedness of the additional 
condition for $\mathscr I\subseteq\mathscr B_T$ to be an ideal of
$\mathscr B_T$.

\begin{thm}\label{thm:GHilb->GGrass}
Let $X=\Spec_S\mathscr B\to S$ be an affine $S$-scheme with 
$G$-operation over $S$. Then the natural morphism of $S$-functors
$\FGHilb_SX\to\FGGrass_S(\mathscr B)$ is a closed embedding. 
\end{thm}
\vspace{-3mm}
\begin{proof}
To show that the canonical inclusion defined above is a closed 
embedding, one has to show that for any $S$-scheme $S'$ and any 
$S'$-valued point of the Grassmannian $[0\to\mathscr H
\stackrel{\phi}{\to}\mathscr B_{S'}\to\mathscr G\to 0]$ 
there is a closed subscheme $Z\subseteq S'$ such that for every 
$S$-morphism $\alpha\colon T\to S'$:
$\alpha$ factors through $Z$ if and only if the quotient 
$[\mathscr B_T\to\alpha^*\mathscr G]\in\FGGrass_S(\mathscr B)(T)$ 
determined by $\alpha$ comes from a quotient of the Hilbert functor.

For $\alpha\colon T\to S'$ the quotient 
$[\alpha^*\mathscr H\stackrel{\alpha^*\phi}{\to}
\mathscr B_{T}\to\alpha^*\mathscr G\to 0]$ of $\O_T$-modules
comes from a quotient of $\mathscr B_T$-modules if and only if the 
image of $\alpha^*\mathscr H$ in $\mathscr B_T$ is a 
$\mathscr B_T$-submodule. Equivalently, the composition 
$\chi_T\colon\mathscr B_T\otimes_{\O_T}\alpha^*\mathscr H\to
\mathscr B_T\otimes_{\O_T}(\alpha^*\phi)(\alpha^*\mathscr H)
\to\mathscr B_T\to\alpha^*\mathscr G$, where the arrow in the middle 
is defined using the multiplication map
$\mathscr B_T\otimes_{\O_T}\mathscr B_T\to\mathscr B_T$ of the 
algebra $\mathscr B_T$, is zero.

We define $Z$, considering the map $\chi_{S'}\colon\mathscr B_{S'}
\otimes_{\O_{S'}}\mathscr H\to\mathscr G$, to be the closed subscheme
with ideal sheaf $\mathscr I$ that is the ideal of $\O_{S'}$ minimal 
with the property $\im(\chi_{S'})\subseteq\mathscr I\mathscr G$.
Locally $\mathscr I$ can be described as follows: 
if $\mathscr G|_U\cong\bigoplus_j\O_U^{(j)}$ for an open $U\subseteq S'$,
then $\mathscr I$ is the ideal sheaf generated by the images of the
coordinate maps $\chi_{U}^{(j)}\colon(\mathscr B_{S'}\otimes_{\O_{S'}}
\mathscr H)|_U\to\O_U^{(j)}$.

This $Z$ has the required property: one has to show that 
$\alpha\colon T\to S'$ factors through $Z$ if and only if
$(\alpha^*\phi)(\alpha^*\mathscr H)\subseteq\mathscr B_T$ is a 
$\mathscr B_T$-submodule.
The map $\chi_T$ can be identified with $\alpha^*\chi_{S'}$. 
The question is local on $S'$, assume that 
$\mathscr G\cong\bigoplus_j\O_{S'}^{(j)}$.
Then there are the equivalences:
\vspace*{-1mm}
\[
\begin{array}{rcl}
(\alpha^*\phi)(\alpha^*\mathscr H)\;\textit{is a $\mathscr B_T$-submodule}
&\Longleftrightarrow& \im(\chi_T)=0\;\textit{in $\alpha^*\mathscr G$}\\
&\Longleftrightarrow& \im(\alpha^*\chi_{S'})=0\;
\textit{in $\alpha^*\mathscr G$}\\
&\Longleftrightarrow& \forall\:j\colon\im(\alpha^*\chi_{S'}^{(j)})=0\;
\textit{in $\O_T^{(j)}$}\\
&\Longleftrightarrow& \forall\:j\colon\im(\chi_{S'}^{(j)})
\subseteq\ker(\O_{S'}\to\alpha_*\O_T)\\
&\Longleftrightarrow&\mathscr I\subseteq\ker(\O_{S'}\to\alpha_*\O_T)\\
&\Longleftrightarrow&\alpha\colon T\to S\;\textit{factors through $Z$}
\end{array}
\]
\vspace*{-8.6mm}

\end{proof}

\begin{cor}\label{cor:reprGHilb,X->Sfinite}
Assume that the Hopf algebra $A$ of $G$ is cosemisimple. 
If $X\to S$ is affine, then the $S$-functor $\FGHilb_SX$ is 
representable. If moreover $X\to S$ is finite, then $\FGHilb_SX$ 
is represented by a projective $S$-scheme.
\end{cor}

\vspace*{-3mm}

\section{Representability of {\boldmath$G$}-Hilbert functors}
\label{sec:reprGHilb}

\vspace*{-2mm}\enlargethispage{1.2mm}
\begin{thm}\label{thm:reprGHilb}
Let $G=\Spec A$ be a finite group scheme over a field $K$ with $A$ 
cosemisimple. Let $X$ be a $G$-scheme algebraic over $K$ and assume 
that a geometric quotient $\pi\colon X\to X/G$, $\pi$ affine, 
of $X$ by $G$ exists.
Then the $G$-Hilbert functor $\FGHilb_KX$ is represented 
by an algebraic $K$-scheme $\GHilb_KX$ and the morphism 
$\tau\colon\GHilb_KX\to X/G$ of proposition \ref{prop:GHilbX->Y} 
is projective.
\end{thm}
\pagebreak

\begin{proof}
Since the group scheme $G$ is finite and $X$ is algebraic, 
$\pi$ is a finite morphism and $X/G$ is algebraic.
Thus, corollary \ref{cor:reprGHilb,X->Sfinite} applies to
the $G$-Hilbert functor $\FGHilb_{X/G}X$ and implies that 
$\FGHilb_{X/G}X$ is represented by an algebraic $K$-scheme
$\GHilb_{X/G}X$ projective over $X/G$.

As explained in remark \ref{rem:GHilboverX/G}, it follows 
that $\FGHilb_KX$ is represented by an algebraic $K$-scheme 
$\GHilb_KX$ and that there is an isomorphism of $K$-schemes 
$\GHilb_KX\cong\GHilb_{X/G}X$ that identifies the morphism 
$\tau\colon\GHilb_KX\to X/G$ with the structure morphism of 
$\GHilb_{X/G}X$, which is a projective morphism.
\end{proof}

\appendix

\section{{\boldmath$G$}-sheaves}\label{sec:Gsheaves}

Let $X$ be a $G$-scheme over $K$ with $s_X\colon G\times_KX\to X$
the group operation. A $G$-equivariant sheaf ($G$-sheaf) on $X$ 
is an $\O_X$-module $\mathscr F$ with the additional structure 
of a $G$-linearisation $\lambda^{\mathscr F}\colon s_X^*\mathscr F
\;\stackrel{\sim}{\longrightarrow}\;p_X^*\mathscr F$ as defined in 
\cite[Ch. 1, \S 3]{MuGIT}.

\begin{rem}\label{rem:Gsheaf}
For a $K$-scheme $T$ and $g\in G(T)=G_T(T)$ there is the diagram

\noindent
\begin{picture}(150,30)
\put(75,24){\makebox(0,0)[c]{$\hspace{14pt} X_T=\,T\times_T X_T=X_T$}}
\put(75,10){\makebox(0,0)[c]{$\hspace{8pt}G_T\times_T X_T$}}
\put(39,4){\makebox(0,0)[c]{$X_T$}}
\put(113,4){\makebox(0,0)[c]{$X_T$}}

\put(75,20){\vector(0,-1){7}}
\put(58,20){\vector(-3,-2){16}}\put(94,20){\vector(3,-2){16}}
\put(63,10){\vector(-4,-1){18}}\put(89,10){\vector(4,-1){18}}

\put(76,17){\makebox(0,0)[l]{\scriptsize$g\times\id_{X_T}$}}
\put(49,17){\makebox(0,0)[r]{\scriptsize$g$}}
\put(103,17){\makebox(0,0)[l]{\scriptsize$\id_{X_T}$}}

\put(55,5){\makebox(0,0)[l]{\scriptsize$s_{X_T}$}}
\put(99,5){\makebox(0,0)[r]{\scriptsize$p_{X_T}$}}
\end{picture}

\noindent
Here the automorphism $g\colon X_T\to X_T$ of $X_T$ over $T$ 
determined by $g\in G(T)$ is denoted as well by $g$.
Pulling back $(\lambda^{\mathscr F})_T\colon s_{X_T}^*\mathscr F_T
\to p_{X_T}^*\mathscr F_T$ from $G_T\times_TX_T$ to $X_T$ by 
$g\times\id_{X_T}$ leads to an isomorphism 
$\lambda^{\mathscr F_T}_g\colon g^*\mathscr F_T\to\mathscr F_T$.
\end{rem}

\begin{prop}\label{prop:Gsheaf}
Let $\mathscr F$ be an $\O_X$-module on a $G$-scheme $X$ over $K$.
Then a $G$-sheaf structure on $\mathscr F$ is equivalent to
the following data: For any $K$-scheme $T$ and any $T$-valued point
$g\in G(T)=G_T(T)$ an isomorphism 
$\lambda^{\mathscr F_T}_g\colon g^*\mathscr F_T\to\mathscr F_T$ 
of $\O_{X_T}$-modules such that 
$\lambda^{\mathscr F_{T'}}_{g_{T'}}=(\lambda^{\mathscr F_T}_g)_{T'}$
for $T$-schemes $T'$ and the following properties are satisfied:\\
\begin{tabular}{rl}
(i)&$\lambda^{\mathscr F_T}_{e_T}=\id_{\mathscr F_T}\;$ for 
$\:e_T\colon T\to G_T$ the identity of $G(T)=G_T(T)$\\
(ii)&$\lambda^{\mathscr F_T}_{hg}=\lambda^{\mathscr F_T}_g\circ
g^*\lambda^{\mathscr F_T}_h\;$ for $\:g,h\in G(T)$\\
\end{tabular}
\end{prop} 

The correspondence is given by the construction in remark
\ref{rem:Gsheaf} and by specialisation to $T=G$, $g=\id_G$, 
that is $\lambda^{\mathscr F}=\lambda^{\mathscr F_G}_{\id_G}$. 
Details can be found in \cite{Bl07}.
The cocycle condition of the original definition arises as 
$\lambda^{\mathscr F_{G\times G}}_{p_1\cdot p_2}=
\lambda^{\mathscr F_{G\times G}}_{p_2}\circ
p_2^*\lambda^{\mathscr F_{G\times G}}_{p_1}$, where 
$p_1,p_2\in G(G\times_K G)$ are the projections.

\pagebreak
\begin{rem}\label{rem:Gsubsheaf}
A subsheaf $\mathscr F'\subseteq\mathscr F$ is called $G$-stable
or a $G$-subsheaf, if \linebreak $\lambda^{\mathscr F}(s_X^*\mathscr F')
\subseteq p_X^*\mathscr F'$ or equivalently 
$\lambda^{\mathscr F_T}_g\colon g^*\mathscr F_T\to\mathscr F_T$ 
restricts to $g^*\mathscr F'_T\to\mathscr F'_T$ for all 
$K$-schemes $T$ and $g\in G(T)$.
If $\mathscr F'\subseteq\mathscr F$ is a $G$-subsheaf, then the
$G$-sheaf structure on $\mathscr F$ restricts to a $G$-sheaf 
structure on $\mathscr F'$. Conversely, an injective $G$-equivariant 
homomorphism determines a $G$-subsheaf.
\end{rem}

We consider the case of trivial $G$-operation on the underlying 
scheme. Comodule structures on sheaves are defined as comodule 
structures on modules, see e.g.\ \cite{Sw}.

\begin{prop}\label{prop:Gsheaf-Acomod}
Let $X$ be a $\,G$-scheme over $K$ with trivial $G$-operation. 
Then for an $\O_X$-module $\mathscr F$ the following data are 
equivalent:\\
\begin{tabular}{rl}
(a)&A $G$-sheaf structure on $\mathscr F$.\\
(b)&An $A$-comodule structure 
$\rho\colon\mathscr F\to A\otimes_K\mathscr F$.\\
\end{tabular}
\end{prop}

The correspondence is given by the adjunction $(p_X^{\;*},p_{X*})$,
for details see \cite{Bl07}. 

\begin{rem}
In particular in the case $X=\Spec L$, $L$ a field, a $G$-sheaf 
structure of a finite dimensional $L$-vector space $V$ is the 
same as an $A$-comodule structure $V\to A\otimes_KV$. This is 
dual to a $\KG$-module structure $\KG\otimes_KV^\vee\to V^\vee$, 
where $\KG=A^\vee$ is the $K$-algebra dual to the coalgebra $A$. 
So we may denote $G$-sheaves on spectra of fields also as
representations.
\end{rem}

\begin{rem}\label{rem:isotypdecomp}
If $A$ is cosemisimple (see e.g. \cite[Ch. XIV]{Sw}), then it 
decomposes into a direct sum $A=\bigoplus_i A_i$ of its simple 
subcoalgebras.
In this case, any $G$-sheaf $\mathscr F$ on a $G$-scheme $X$ 
with trivial $G$-operation decomposes into a direct sum 
$\mathscr F=\bigoplus_i\mathscr F_i$ of isotypic components 
$\mathscr F_i$, whose comodule structure reduces to an 
$A_i$-comodule structure (the $\mathscr F_i$ are given as 
preimages $\rho^{-1}(A_i\otimes_K\mathscr F)$).
For points $x\in X$ this decomposition induces the isotypic
decomposition of the representations on the fibres $\mathscr F_x$.
\end{rem}

\bigskip
\noindent
Mark Blume\\
Mathematisches Institut, Universit\"at M\"unster,\\
Einsteinstrasse 62, 48149 M\"unster, Germany\\
E-mail: mark.blume@uni-muenster.de

\end{document}